\documentclass[12pt,a4paper]{amsart}
\usepackage{enumerate}
\usepackage{amsmath,amsthm,verbatim,amssymb,amsfonts,amscd,graphicx}
\usepackage{graphics}
\usepackage{hyperref}
\usepackage{mathrsfs}
\usepackage{tikz}
\usepackage{tikz-cd}
\usepackage[margin=2.9cm]{geometry}
\usepackage[latin2]{inputenc}
\usepackage{xcolor}
\colorlet{lightCyan}{cyan!35}
\colorlet{deepCyan}{cyan!75}

\newcommand{\X}{\mathbf{X}}
\newcommand{\ii}{\mathbf{i}}
\newcommand{\nn}{\mathbf{n}}

\newcommand{\mm}{\mathbf{m}}
\newcommand{\II}{\mathbf{I}}
\newcommand{\JJ}{\mathbf{J}}

\newcommand{\ZZ}{\mathbb{Z}}
\newcommand{\NN}{\mathbb{N}}
\newcommand{\R}{\mathcal{R}}

\theoremstyle{plain}
\newtheorem{theorem}{Theorem}[section]

\newtheorem{corollary}{Corollary}[section]
\newtheorem*{corollary*}{Corollary}
\newtheorem{lemma}{Lemma}[section]
\newtheorem{proposition}{Proposition}[section]
 
\newtheorem{question}{Question}[section] 
\newtheorem{example}{Example}[section]
\newtheorem*{example*}{Example}

\newtheorem{remarks}{Remarks}[section]
\theoremstyle{definition}

\newtheorem*{definition*}{Definition}
\newcommand{\thickfrac}[2]{\genfrac{}{}{0.75pt}{}{#1}{#2}}

\begin{document}
\title{Finitely generated saturated multi-Rees algebras}
\author{Suprajo Das}
\address{Indian Institute of Technology Bombay, Mumbai, 400076, India}
\email{dassuprajo@gmail.com}
\author{Sudeshna Roy}
\address{Tata Institute of Fundamental Research, Mumbai, 400005, India}
\email{sudeshnaroy.11@gmail.com}

\begin{abstract}
We study the question of finite generation of saturated multi-Rees algebras and investigate the asymptotic behaviour of related length functions. In the setup of excellent local domains, we show that the saturated multi-Rees algebra of a finite collection of ideals is finitely generated when the analytic spread is not maximal and the associated length function eventually agrees with a polynomial. Similar results are obtained when we restrict to two-dimensional local UFDs with no restrictions on the analytic spread. We further prove that the saturated multi-Rees algebra of finitely many monomial ideals in a polynomial ring modulo an irreducible monomial ideal, is always finitely generated. In this case, the corresponding length function is shown to exhibit piecewise quasi-polynomial behaviour. We also produce multi-ideal versions of a theorem of Amao.
\end{abstract}

\maketitle

\section{Introduction}
The inspiration of this manuscript mainly comes from three sources: \cite{DC4} by Cutkosky, Herzog and Srinivasan, \cite{J} by Herzog, Puthenpurakal and Verma and \cite{HHT} by Herzog, Hibi and Trung. Here we produce multi-ideal versions of some of the results obtained in the aforementioned papers.

Suppose that $(R,m_R)$ is a Noetherian local ring (or a standard graded Noetherian ring with graded maximal ideal) of dimension $d$ and $I_1,\ldots,I_r$ are ideals in $R$. We are interested in studying the asymptotic behaviour of the function $l_R\left(H^0_{m_R}\left(R/I_1^{n_1}\cdots I_r^{n_r}\right)\right)$, where $n_1,\ldots,n_r\in \NN$. In Theorem \ref{rightbound} we show that this length function is bounded above by a numerical polynomial in $n_1,\ldots,n_r$ of total degree at most $d$. If we only consider a single ideal $I$ then it is known previously that
\begin{equation}\label{limitornot}
\varepsilon(I):=\limsup\limits_{n\to\infty} \dfrac{l_R\left(H^0_{m_R}\left(R/I^n\right)\right)}{n^d/d!}<\infty,
\end{equation}
which is a special case of our result. Katz and Validashti \cite{KV} calls this invariant $\varepsilon(I)$ the \emph{epsilon multiplicity} of the ideal $I$. Now it begs the question if it is possible to define `\emph{mixed}' versions of the epsilon multiplicity.

The associated length function $l_R\left(H^0_{m_R}\left(R/I^n\right)\right)$ can exhibit very mysterious behaviour. An example constructed in \cite{DC3} shows that the invariant $\varepsilon(I)$ can be an irrational number. In general, the sequence $\left\{l_R\left(H^0_{m_R}\left(R/I^n\right)\right)\right\}_{n\in\mathbb{N}}$ does not even have quasi-polynomial growth eventually and hence cannot arise from a finitely generated graded module. A very interesting question is to decide whether the $\limsup$ in \eqref{limitornot} can be replaced by a limit. Cutkosky \cite{DC2} has given a positive answer to this question when $R$ is analytically unramified. A related problem is to describe instances where this limit can be a rational number. It turns out that if the saturated Rees algebra $$S_{m_R}(I):=\bigoplus_{n\in\mathbb{N}}(I^n\colon_R m_R^{\infty})$$ is finitely generated then $\varepsilon(I)$ is rational. However, the saturated Rees algebra need not always be finitely generated. More generally, we consider the $\mathbb{N}^r$-graded $R$-algebra $$S_J(I_1,\ldots,I_r) = \bigoplus\limits_{(n_1,\ldots,n_r)\in\mathbb{N}^r}\left(I_1^{n_1}\cdots I_r^{n_r}\colon_R J^{\infty}\right)$$ where $I_1,\ldots,I_r, J$ are ideals in $R$. We say that $S_J(I_1,\ldots,I_r)$ is the \emph{saturated multi-Rees algebra} of $I_1,\ldots,I_r$ with respect to $J$. If $J=m_R$ then $S_{m_R}(I_1,\ldots,I_r)$ will be called the \emph{saturated multi-Rees algebra} of $I_1,\ldots,I_r$.

The fundamental difficulty in working with multigraded Noetherian algebras is that it might not contain a standard multigraded Veronese subalgebra which is also Noetherian. It is worth mentioning that the results in Sections \ref{amaogen}, \ref{prooffirst} and \ref{dimensiontwo} of this manuscript have a `\emph{diagonal}' flavour. We now describe certain situations when $S_J(I_1,\ldots,I_r)$ is a finitely generated $\NN^r$-graded $R$-algebra.
  
\begin{theorem}\label{mainthm1}
Suppose that $(R,m_R)$ is an excellent local domain of dimension $d$ and $I_1,\ldots,I_r,J$ are ideals in $R$. Assume that the analytic spread $\ell({I_1}_P\cdots {I_r}_P) < \dim R_P$ for all $P\in V(J)$. Then the saturated multi-Rees algebra $S_J(I_1,\ldots,I_r)$ of $I_1,\ldots,I_r$ with respect to $J$ is Noetherian. There exists a numerical polynomial $F \in \mathbb{Q}[n_1,\ldots,n_r]$ of total degree at most $\ell(I_1\cdots I_r) + \dim R/J - 1$ such that $$e\left(H^0_{J}\left(R/I_1^{n_1}\cdots I_r^{n_r}\right)\right)=F(n_1,\ldots,n_r)$$ for all $n_1,\ldots,n_r\gg0$.
\end{theorem}

The first statement of Theorem \ref{mainthm1} was already proven for $r=1$ by Cutkosky, Herzog and Srinivasan in \cite{DC4}. By generalizing their proof we show that $S_J(I_1,\ldots,I_r)$ is contained inside the integral closure of the multi-Rees algebra $\mathcal{R}(I_1,\ldots,I_r)$, thereby proving our first assertion. The polynomial behaviour is a consequence of the fact that $S_{J}(I_1,\ldots,I_r)$ is a finite module over $\mathcal{R}(I_1,\ldots,I_r)$, see Theorem \ref{finitely} for the precise statements.

An immediate consequence of Theorem \ref{mainthm1} shows that if we consider ideals $I_1,\ldots,I_r$ in an excellent local domain $(R,m_R)$ such that the analytic spread $\ell(I_1\cdots I_r)$ is not maximal, then the saturated multi-Rees algebra $S_{m_R}(I_1,\ldots,I_r)$ is Noetherian. The corresponding length function $l_R\left(H^0_{m_R}\left(R/I_1^{n_1}\cdots I_r^{n_r}\right)\right)$ is eventually given by a numerical polynomial of total degree at most $\ell(I_1\cdots I_r)-1$, see Corollary \ref{smallpolynomial} for details. The next result concerns the saturated multi-Rees algebra of monomial ideals in a polynomial ring modulo an irreducible monomial ideal.

\begin{theorem}\label{mainthm2}
Suppose that $R = K[X_1,\ldots,X_d]$ is a standard graded polynomial ring in $d$ variables over a field $K$. Consider the graded ring $S=R/L$ where $L\subseteq R$ is either the zero ideal or an irreducible monomial ideal. Suppose that $I_1,\ldots,I_r$ are ideals generated by monomials in $S$ and $J$ is any ideal in $S$. Then the saturated multi-Rees algebra $S_J(I_1,\ldots,I_r)$ of $I_1,\ldots,I_r$ with respect to $J$ is Noetherian. If $J=m_S$, where $m_S$ is the graded maximal ideal of $S$, then there exists a piecewise quasi-polynomial $F \colon \mathbb{N}^r \to \mathbb{N}$ such that $$l_S\left(H^0_{m_S}\left(S/I_1^{n_1}\cdots I_r^{n_r}\right)\right) = F(n_1,\ldots,n_r)$$ for all $n_1,\ldots,n_r\gg0$. Moreover, this length function is bounded above by a numerical polynomial in $n_1,\ldots,n_r$ of total degree at most $\dim S$ for all $n_1,\ldots,n_r\gg0$.
\end{theorem}

We shall prove this result in Section \ref{proofsecond}. In the case $r=1$, $L=(0)$ and $J\subset R$ a monomial ideal, Herzog, Hibi and Trung \cite{HHT} have shown that the saturated Rees algebra $S_{J}(I)$ is Noetherian and Herzog, Puthenpurakal and Verma \cite{J} proved that the corresponding function $e\left(H^0_{J}\left(R/I^n\right)\right)$ is eventually a quasi-polynomial with constant leading term. Theorem \ref{mainthm2} and Corollary \ref{maincor} generalies these existing results but our methods differ from \cite{HHT} and \cite{J}. We use a projection argument to obtain the first assertion of our theorem. For proving the second assertion we invoke a result from \cite{das2} obtained by the first author using Presburger arithmetic. We provide some examples to support our results.

In section \ref{dimension2}, we concern ourselves with two-dimensional local domains. More precisely, we show that if the ambient ring $R$ is a local UFD of dimension two then the saturated multi-Rees algebra of finitely many ideals in $R$ is Noetherian and the associated length function is given by a polynomial.

\section{Notations and Preliminaries}
We begin by fixing some notations and by recalling some basic facts about multi-Rees rings.

Let $\NN$ be the set of natural numbers including zero. Vectors in $\ZZ^r$ will be written by bold-faced letters, e.g., $\mathbf{a},\nn,\mathbf{i},\mathbf{j},\mathbf{h},\ldots$. The norm of a multi-index $\nn = (n_1,\ldots,n_r) \in\ZZ^r$ is $|\nn| = n_1+\cdots +n_r$. We denote by $\mathbf{0}$ and $\mathbf{1}$ the vectors $(0,\ldots,0)\in\NN^r$ and $(1,\ldots,1)\in\NN^r$ respectively. For every $i=1,\ldots,r$ let $\mathbf{e}_i$ denote the standard unit vector in $\NN^r$ with $1$ at the $i$-th spot and zeros elsewhere. If $m_i > n_i$ (resp. $m_i\geq n_i$) for all $i=1,\ldots,r$ then we write $\mm > \nn$ (resp. $\mm \geq \nn$). If $n_i\gg 0$ for all $i=1,\ldots,r$ then we write $\nn \gg \mathbf{0}$. For $\ii,\nn\in\NN^r$ we denote by $\nn^{\ii} = n_1^{i_1}\cdots n_r^{i_r}$.

Let $S = \bigoplus_{\nn\in\NN^r}S_{\nn}$ be an $\NN^r$-graded ring. If $\nn\in\NN^r$ and $s\in S_{\nn}$, we say that $s$ is \emph{homogeneous} of degree $\nn$. We say that $S$ is a \emph{standard} $\NN^r$-graded $S_{\mathbf{0}}$-algebra if $S$ is generated over $S_{\mathbf{0}}$ by elements in $S_{\mathbf{e}_1},\ldots,S_{\mathbf{e}_r}$.

Suppose that $R$ is a ring, $I_1,\ldots,I_r$ are ideals in $R$ and $t_1,\ldots,t_r$ are variables over $R$. Given $\nn\in\NN^r$ we write $\mathbf{I}^{\nn} = I_1^{n_1}\cdots I_r^{n_r}$ and $\mathbf{t}^{\nn} = t_1^{n_1}\cdots t_r^{n_r}$ for short. The \emph{multi-Rees algebra of $I_1, \ldots, I_r$} is defined as $$\R(I_1,\ldots,I_r) := \bigoplus\limits_{\nn\in\mathbb{N}^r}\mathbf{I}^{\nn}\mathbf{t}^{\nn}$$ and $\overline{\R(I_1,\ldots,I_r)}$ denotes the integral closure of $\R(I_1,\ldots,I_r)$ in its total ring of fractions. If we further assume that $R$ is a Noetherian local ring with unique maximal ideal $m_R$ then the analytic spread $\ell(I_1,\ldots,I_r)$ of $I_1,\ldots,I_r$ is defined as $$\ell(I_1,\ldots,I_r):= \dim \left(\dfrac{\R(I_1,\ldots,I_r)}{m_R\R(I_1,\ldots,I_r)}\right) = \dim \left(\bigoplus\limits_{\nn\in\NN^r}\dfrac{\II^{\nn}}{m_R\II^{\nn}}\right).$$ Given an ideal $I\subseteq R$, let $\overline{I}$ denote the integral closure of $I$ in $R$.

\section{Multigraded Rees algebras and Hilbert functions}
The following result is well-known and describes the integral closure of multi-Rees algebras. This is a multi-ideal version of \cite[Proposition 5.2.1]{HS} and \cite[Proposition 5.2.4]{HS}.

\begin{lemma}\label{reesint}
 Suppose that $R$ is a ring, $I_1,\ldots,I_r$ are ideals in $R$ and $t_1,\ldots,t_r$ are variables over $R$. Then the integral closure of the multi-Rees algebra $\R(I_1,\ldots,I_r)$ in $R[t_1,\ldots,t_r]$ equals the $\mathbb{N}^r$-graded ring $$\bigoplus\limits_{\nn\in\mathbb{N}^r} \overline{\II^{\nn}}\mathbf{t}^{\nn}.$$ Let $\overline{R}$ denote the integral closure of $R$ in its total ring of fractions. Then the integral closure of $\R(I_1,\ldots,I_r)$ in its total ring of fractions equals the $\mathbb{N}^r$-graded ring $$\bigoplus\limits_{\nn\in\mathbb{N}^r} \overline{\II^{\nn}\overline{R}}\mathbf{t}^{\nn}.$$
\end{lemma}
\begin{proof}
  Let $S$ denote the integral closure of $\R(I_1,\ldots,I_r)$ in $R[t_1,\ldots,t_r]$. By \cite[Theorem 2.3.2]{HS} $S$ is an $\NN^r$-graded submodule of $R[t_1,\ldots,t_r]$. Let $s\in S_{\nn}$ where $\nn\in\NN^r$. Write $s = s_{\nn}\mathbf{t}^{\nn}$ for some $s_{\nn} \in R$. As $s$ is integral over $\R(I_1,\ldots,I_r)$, there exists an integer $e>0$ and $a_i \in \R(I_1,\ldots,I_r)$ for all $i=1,\ldots,e,$ such that $$s_{\nn}^e \mathbf{t}^{e\cdot\nn} + a_1s_{\nn}^{e-1}\mathbf{t}^{(e-1)\cdot\nn} + \cdots + a_e = 0.$$ Expand each $$a_i = \sum_{\mathbf{j}\in\mathbb{N}^r}a_{i,\mathbf{j}}\mathbf{t}^{\mathbf{j}}$$ as a finite sum with $a_{i,\mathbf{j}} \in \II^{\mathbf{j}}$. The homogeneous part of degree $e\cdot\nn$ in the equation above is exactly $$\mathbf{t}^{e\cdot\nn}\left(s_{\nn}^e + a_{1,\nn}s_{\nn}^{e-1} + \cdots + a_{e,e\cdot\nn}\right) = 0.$$ As $a_{i,i\cdot\nn} \in \left(\II^{\nn}\right)^i$, this equation says that $s_{\nn}$ is integral over the ideal $\II^{\nn}$. The other inclusion is easy to prove. This proves the first statement of the lemma.
  
  For the proof of the second statement, observe that for all $\nn\in\NN^r$, we have $\overline{\II^{\nn}\overline{R}} = \overline{\overline{\II^{\nn}}\;\overline{R}}$. The integral closure of $\R(I_1,\ldots,I_r)$ in its total ring of fractions clearly contains $\overline{R}[I_1\overline{R}t_1,\ldots,I_r\overline{R}t_r]$. By the first statement of the lemma, we know that integral closure of $\overline{R}[I_1\overline{R}t_1,\ldots,I_r\overline{R}t_r]$ in $\overline{R}[t_1,\ldots,t_r]$ is $$\bigoplus\limits_{\nn\in\mathbb{N}^r} \overline{\II^{\nn}\overline{R}}\mathbf{t}^{\nn}.$$ But $\overline{R}[t_1,\ldots,t_r]$ is integrally closed in its total ring of fractions, so that the integral closure of $\R(I_1,\ldots,I_r)$ is as displayed above.
\end{proof}

\begin{proposition}\cite[Lemma 3.4]{ghosh}\label{constant}
Suppose that $(R,m_R)$ is a Noetherian local ring, $A=\bigoplus_{\mathbf{n}\in\mathbb{N}^r} A_{\mathbf{n}}$ is a standard $\mathbb{N}^r$-graded Noetherian $R$-algebra with $A_{\mathbf{0}}=R$ and $M=\bigoplus_{\mathbf{n}\in\mathbb{N}^r} M_{\mathbf{n}}$ is a finitely generated $\mathbb{N}^r$-graded $A$-module. There exists a vector $\mathbf{h}\in\mathbb{N}^r$ such that $\mathrm{Ann}_R M_{\mathbf{n}} = \mathrm{Ann}_R M_{\mathbf{h}}$ for all $\mathbf{n}\geq \mathbf{h}$. In particular, $\dim M_{\mathbf{n}} = \dim M_{\mathbf{h}}$ for all $\mathbf{n}\geq \mathbf{h}$.
\end{proposition}

In the setup of Proposition \ref{constant}, $t:=\dim_R M_{\mathbf{h}}$ is called the \emph{limit dimension} of the $R$-modules $M_{\nn}$. The Hilbert-Samuel multiplicity of $M_{\nn}$ is given by $$e(M_{\nn}):=\lim\limits_{k\to\infty}\dfrac{l_R\left(M_{\nn}/m_R^k M_{\nn}\right)}{k^t/t!}.$$ If $t=0$ then the modules $M_{\nn}$ eventually become finite length $R$-modules. In that case we have $e(M_{\nn})=l_R(M_{\nn})$.

The following result generalises an observation due to Katz and Validashti \cite[Observation 4.1]{KV} for multiple ideals.

\begin{theorem}\label{rightbound}
Suppose that $(R,m_R)$ is a Noetherian local ring of dimension $d$ and $I_1,\ldots,I_r$ are ideals in $R$. Then there exists a numerical polynomial $F\in\mathbb{Q}[n_1,\ldots,n_r]$ of total degree at most $d$ such that $$l_R\left(H^0_{m_R}\left(R/\II^{\nn}\right)\right) \leq F(\nn)$$ for all $\nn\gg\mathbf{0}$.
\end{theorem}

\begin{proof}
Given an $r$-tuple $\nn \in \mathbb{N}^r$ and an integer $k\in\{1,\ldots,r\}$, we define $$\mathbb{I}(\nn, k):= \prod_{i=k}^r I_i^{n_i}.$$ Also set $\mathbb{I}(\nn, r+1):= R$ and notice that $\mathbb{I}(\nn, k)= I_k^{n_k} \mathbb{I}(\nn, k+1)$ for all $1\leq k\leq r$. Consider the chain of inclusions of the ideals
\[\mathbb{I}(\nn, 1) \subseteq I_1^{n_1-1}\mathbb{I}(\nn, 2)\subseteq \cdots \subseteq \mathbb{I}(\nn, 2) \subseteq I_2^{n_2-1}\mathbb{I}(\nn, 3)\subseteq \cdots \subseteq \mathbb{I}(\nn, r) \subseteq \cdots \subseteq I_r \subseteq R.\]
The subadditivity of $l_R\left(H^0_{m_R}\left(-\right)\right)$ on short exact sequences shows that
\begin{equation}\label{rel_jmul_emul}
l_R\left(H^0_{m_R}\left(R/\II^{\nn}\right)\right) \leq \sum_{k=1}^{r}\sum_{t=0}^{n_k-1} l_R\left(H^0_{m_R}\left(\dfrac{I_k^t\mathbb{I}(\nn, k+1)}{I_k^{t+1}\mathbb{I}(\nn, k+1)}\right)\right).
\end{equation}
Consider the $\mathbb{N}^{r-k+1}$-graded multi-Rees algebra
\[\mathcal{R}_k:=\bigoplus_{\mm_k:=(m_k,\ldots,m_r) \in \mathbb{N}^{r-k+1}} I_{k}^{m_k}\cdots I_r^{m_r} \] of the ideals $I_k,\ldots,I_r$ and the associated graded ring $\mathcal{G}_k:=\mathcal{R}_k/I_k \mathcal{R}_k$. Then $H^0_{m_R} \left(\mathcal{G}_k\right)$ is a finitely generated $\NN^{r-k+1}$-graded $\mathcal{G}_k$-module which is annihilated by $m_R^{c}$ for some integer $c>0$. Therefore, $H^0_{m_R} \left(\mathcal{G}_k\right)$ can be regarded as a finitely generated $\NN^{r-k+1}$-graded module over $\mathcal{G}_k/m_R^c \mathcal{G}_k$, which is a standard $\NN^{r-k+1}$-graded Noetherian algebra over the Artinian local ring $R/(m_R^c+I_k)$. It then follows from \cite[Lemma $1.1$ and Theorem $4.1$]{HHRT} that there exists a numerical polynomial $P_{k}\in \mathbb{Q}[m_k,\ldots,m_r]$ of total degree at most $\dim \mathcal{G}_k - (r-k+1) \leq d-1$ such that $$l_R\left(\left(H^0_{m_R}\left(\mathcal{G}_k\right)\right)_{\mm_k}\right) = l_R\left(H^0_{m_R}\left(\dfrac{I_k^{m_k}\cdots I_r^{m_r}}{I_k^{m_k+1}I_{k+1}^{m_{k+1}}\cdots I_r^{m_r}}\right)\right) = P_{k}\big(m_{k},\ldots,m_{r}\big)$$ for all $m_k,\ldots,m_r \gg 0$. Consider the sum $$F_k(n_k,\ldots,n_r) := \sum\limits_{t=0}^{n_k-1} l_R\left(H^0_{m_R}\left(\dfrac{I_k^t\mathbb{I}(\nn, k+1)}{I_k^{t+1}\mathbb{I}(\nn, k+1)}\right)\right)$$ and observe that the difference function \[F_k(n_k,\ldots,n_r) - F_k(n_k-1,n_{k+1},\ldots,n_r) = l_R\left(H^0_{m_R}\left(\dfrac{I_k^{n_k-1}\mathbb{I}(\nn, k+1)}{I_k^{n_k}\mathbb{I}(\nn, k+1)}\right)\right)\] eventually agrees with the numerical polynomial $P_k(n_k,\cdots,n_r)$ having total degree at most $d-1$. We now apply \cite[Proposition $2.3.6$]{nim} to conclude that each $F_k(n_k,\ldots,n_r)$ is eventually given by a numerical polynomial in $n_k,\ldots,n_r$ of total degree at most $d$. Finally we use \eqref{rel_jmul_emul} to reach our desired conclusion.
\end{proof}

\section{Multigraded generalisations of a theorem of Amao}\label{amaogen}

The following theorem extends Amao's result \cite[Theorem 3.2]{JO} to a finite number of ideals.

\begin{theorem}\label{rarepolynomial}
Suppose that $(R,m_R)$ is a Noetherian local ring of dimension $d$ and $I_1\subseteq J_1, \ldots, I_r\subseteq J_r$ are ideals in $R$ such that $l_R(J_i/I_i)<\infty$ for all $i=1,\ldots,r$. Then the multigraded length function $l_R\left(\mathbf{J}^{\nn}/\II^{\nn}\right)$ eventually agrees with a numerical polynomial in $\mathbb{Q}[n_1,\ldots,n_r]$ of total degree at most $d$ and all monomials of highest degree in this polynomial have nonnegative coefficients. In particular, $$\mathrm{total\; degree}\; l_R\left(\JJ^{\nn}/\II^{\nn}\right) = \mathrm{deg}\; l_R\left(\left(J_1\cdots J_r\right)^n/\left(I_1\cdots I_r\right)^n\right)\leq d.$$
\end{theorem}
\begin{proof}
Given an $r$-tuple $\nn \in \mathbb{N}^r$ and an integer $k \in \{1, \ldots, r\}$, we define \[\mathbb{J}(k,\nn):= \prod_{j=1}^{k}J_j^{n_j} \quad \mbox{ and } \quad \mathbb{I}(\nn, k):= \prod_{i=k}^{r}I_i^{n_i}.\] Also set $\mathbb{J}(0,\nn) := R$ and $\mathbb{I}(\nn,r+1) := R$. There exist inclusion of ideals $$\mathbb{J}(k-1,\nn) J_k^{t-1} I_k^{n_k-t+1} \mathbb{I}(\nn, k+1)  \subseteq \mathbb{J}(k-1,\nn) J_k^t I_k^{n_k-t} \mathbb{I}(\nn, k+1)$$ for each $k=1,\ldots,r$ and for all $t=1,\ldots,n_k$. Using the additivity of lengths on short exact sequences, we can write
\begin{equation}\label{amaosum}
 l_R\left(\mathbf{J}^{\nn}/\II^{\nn}\right) = \sum\limits_{k=1}^r \sum\limits_{t=1}^{n_k}l_R\left(\dfrac{\mathbb{J}(k-1,\nn) J_k^t I_k^{n_k-t} \mathbb{I}(\nn, k+1)}{\mathbb{J}(k-1,\nn) J_k^{t-1} I_k^{n_k-t+1} \mathbb{I}(\nn, k+1)}\right).
\end{equation}
Given $\nn \in \mathbb{N}^r, (u,v) \in \mathbb{N}^2$ and $k \in \{1, \ldots, r\}$, we define an $(r+1)$-tuple by \[\left(\nn(k),u,v\right):=(n_1,\ldots,n_{k-1},u,v,n_{k+1},\ldots,n_r) \in \mathbb{N}^{r+1}.\] For each $k\in\{1,\ldots,r\}$ we consider the $\mathbb{N}^{r+1}$-graded multi-Rees algebra
$$\mathcal{R}_k := \bigoplus_{\left(\nn(k),u,v\right)\in\mathbb{N}^{r+1}}\mathbb{J}(\nn, k-1) J_k^u I_k^{v} \mathbb{I}(\nn, k+1)$$ of $r+1$ ideals $J_1,\ldots,J_{k},I_k,\ldots,I_r$. Also consider the $\mathbb{N}^{r+1}$-graded ideals
	\begin{align*}
	I_k\mathcal{R}_k &= \bigoplus\limits_{\substack{\left(\nn(k),u,v\right)\in\mathbb{N}^{r+1}\\ u\geq 1}}\mathbb{J}(k-1,\nn) J_k^{u-1} I_k^{v+1} \mathbb{I}(\nn, k+1),\\
	J_k\mathcal{R}_k &=  \bigoplus\limits_{\substack{\left(\nn(k),u,v\right)\in\mathbb{N}^{r+1},\\ u\geq 1}}\mathbb{J}(k-1,\nn) J_k^{u} I_k^{v} \mathbb{I}(\nn, k+1)
	\end{align*}
	of $\mathcal{R}_k$ and their quotient $$J_k\mathcal{R}_k/I_k\mathcal{R}_k = \bigoplus\limits_{\substack{\left(\nn(k),u,v\right)\in\mathbb{N}^{r+1},\\ u\geq 1}} \dfrac{\mathbb{J}(k-1,\nn)J_k^{u} I_k^{v} \mathbb{I}(\nn, k+1)}{\mathbb{J}(k-1,\nn) J_k^{u-1} I_k^{v+1} \mathbb{I}(\nn, k+1)}.$$ We now impose an $\mathbb{N}^r$-grading on $\mathcal{R}_k$ by setting $$\left(\mathcal{R}_k\right)_{\nn}:= \bigoplus\limits_{u+v=n_k}\mathbb{J}(k-1,\nn) J_k^u I_k^{v} \mathbb{I}(\nn, k+1)$$ and $\mathcal{R}_k$ becomes a standard $\NN^r$-graded Noetherian $R$-algebra with this new grading. Observe that $J_k\mathcal{R}_k/I_k\mathcal{R}_k$ becomes a finitely generated $\mathbb{N}^{r}$-graded $\mathcal{R}_k$-module, where the $\nn$-th component of $J_k\mathcal{R}_k/I_k\mathcal{R}_k$ is given by $$\left(J_k\mathcal{R}_k/I_k\mathcal{R}_k\right)_{\nn}:= \bigoplus\limits_{\substack{u+v=n_k,\\u\geq 1}} \dfrac{\mathbb{J}(k-1,\nn)J_k^{u} I_k^{v} \mathbb{I}(\nn, k+1)}{\mathbb{J}(k-1,\nn)J_k^{u-1} I_k^{v+1} \mathbb{I}(\nn, k+1)}.$$ Moreover, each component of $J_k\mathcal{R}_k/I_k\mathcal{R}_k$ is an $R$-module of finite length. Therefore $J_k\mathcal{R}_k/I_k\mathcal{R}_k$ is annihilated by $m_R^{c}$ for some integer $c>0$. So $J_k\mathcal{R}_k/I_k\mathcal{R}_k$ can be regarded as a finitely generated $\mathbb{N}^{r}$-graded module over $\mathcal{R}_k/m_R^c\mathcal{R}_k$ and the latter is a standard $\mathbb{N}^{r}$-graded Noetherian algebra over the Artinian local ring $R/m_R^c$. From \cite[Theorem $4.1$]{HHRT} we know that for each $k\in\{1,\ldots,r\}$, there exists a numerical polynomial $F_k \in \mathbb{Q}[n_1,\ldots,n_r]$ such that $$l_R\left(\left(\dfrac{J_k\mathcal{R}_k}{I_k\mathcal{R}_k}\right)_{\nn}\right) = \sum\limits_{u=1}^{n_k} l_R\left(\dfrac{\mathbb{J}(k-1,\nn)J_k^{u} I_k^{n_k-u} \mathbb{I}(\nn, k+1)}{\mathbb{J}(k-1,\nn)J_k^{u-1} I_k^{n_k-u+1} \mathbb{I}(\nn, k+1)}\right) = F_k(\nn)$$ for all $\nn \gg \mathbf{0}$ and all monomials of highest degree in this polynomial have nonnegative coefficients. From \eqref{amaosum} it now follows that our desired length function $l_R\left(\mathbf{J}^{\nn}/\II^{\nn}\right)$ agrees with the numerical polynomial $F(\nn):=\sum_{k=1}^r F_k(\nn)$ for all $\nn \gg \mathbf{0}$ and all monomials of highest degree in $F(\nn)$ also have nonnegative coefficients. From Amao's theorem \cite[Theorem 3.2]{JO} we know that there exists a numerical polynomial $G(n)\in \mathbb{Q}[n]$ of degree at most $d$ such that $$l_R\left(\left(J_1\cdots J_r\right)^n/\left(I_1\cdots I_r\right)^n\right) = G(n)$$ for all $n\gg 0$. The last assertion then follows from the observation that $F(n,\ldots,n) = G(n)$ for all $n\gg 0$.
\end{proof}

The next lemma ensures existence of limit dimensions of certain modules, also see \cite[Lemma $2.2$]{cat}.

\begin{lemma}\label{eventualdim}
Suppose that $R$ is a Noetherian ring and $I_1\subseteq J_1, \ldots, I_r\subseteq J_r$ are ideals in $R$. Then there exists an $r$-tuple $\mathbf{h}\in\NN^r$ such that $$\sqrt{\left(\II^{\nn}\colon_R \JJ^{\nn}\right)} = \sqrt{\left(\II^{\mathbf{h}}\colon_R \JJ^{\mathbf{h}}\right)}$$ for all $r$-tuples $\nn \geq \mathbf{h}$.
\end{lemma}
\begin{proof}
Let $x\in \sqrt{\left(\II^{\nn}\colon_R \JJ^{\nn}\right)}$ then $x^c \JJ^{\nn} \subset \II^{\nn}$ for some integer $c>0$. Therefore $$x^{2c} \JJ^{\nn + \mathbf{e}_i} = x^c x^c \JJ^{\nn}J_i \subseteq x^c\II^{\nn}J_i \subseteq x^c \JJ^{\nn}I_i \subseteq \II^{\nn}I_i = \II^{\nn + \mathbf{e}_i}$$ for all $i=1,\ldots,r$. In other words, $\sqrt{\left(\II^{\nn}\colon_R \JJ^{\nn}\right)} \subseteq \sqrt{\left(\II^{\nn^{\prime}}\colon_R \JJ^{\nn^{\prime}}\right)}$ whenever $\nn^{\prime}\geq \nn$. Now consider the set of ideals $$\Lambda := \left\{\sqrt{\left(\II^{\nn}\colon_R \JJ^{\nn}\right)} \mid \nn\in\NN^r\right\}$$ which is partially ordered by inclusion. Since $R$ is Noetherian, the set $\Lambda$ has a maximal element with respect to inclusion. It can be easily verified that the maximal element of $\Lambda$ is unique as $\sqrt{\left(\II^{\nn}\colon_R \JJ^{\nn}\right)} \subseteq \sqrt{\left(\II^{\nn^{\prime}}\colon_R \JJ^{\nn^{\prime}}\right)}$ whenever $\nn^{\prime}\geq \nn$. Let $\sqrt{\left(\II^{\mathbf{h}}\colon_R \JJ^{\mathbf{h}}\right)}$ be the maximum element of $\Lambda$ and our result follows.
\end{proof}

The following result further generalises Theorem \ref{rarepolynomial} by producing a multigraded version of \cite[Proposition $2.4$]{cat} or \cite[Proposition $2.1$]{J}.

\begin{theorem}
Suppose that $(R,m_R)$ is a Noetherian local ring of dimension $d$ and $I_1\subseteq J_1, \ldots, I_r\subseteq J_r$ are ideals in $R$. In view of Lemma \ref{eventualdim}, set $$t := \lim\limits_{\nn\to\infty} \dim \dfrac{R}{\sqrt{\II^{\nn}\colon_R \JJ^{\nn}}} = \lim\limits_{n\to\infty}\dim \dfrac{R}{\sqrt{(I_1\cdots I_r)^n\colon_R (J_1\cdots J_r)^{n}}}.$$ Then the multigraded multiplicity function $e\left(\mathbf{J}^{\nn}/\II^{\nn}\right)$ eventually agrees with a numerical polynomial in $\mathbb{Q}[n_1,\ldots,n_r]$ of total degree at most $d-t$ and all monomials of highest degree in this polynomial have nonnegative coefficients. In particular, $$\mathrm{total\; degree}\; e\left(\JJ^{\nn}/\II^{\nn}\right) = \mathrm{deg}\; e\left(\left(J_1\cdots J_r\right)^n/\left(I_1\cdots I_r\right)^n\right)\leq d-t.$$
\end{theorem}
\begin{proof}
From Lemma \ref{eventualdim} we know that there exists an $r$-tuple $\mathbf{h}\in\NN^r$ such that $\sqrt{\left(\II^{\nn}\colon_R \JJ^{\nn}\right)} = \sqrt{\left(\II^{\mathbf{h}}\colon_R \JJ^{\mathbf{h}}\right)}$ for all $\nn \geq \mathbf{h}$. Then $t = \dim \left(\II^{\nn}/\JJ^{\nn}\right)$ for all $\nn \geq \mathbf{h}$. Using the associativity formula, for $\nn \geq \mathbf{h}$ we have
\begin{equation*}
e\left(\mathbf{J}^{\nn}/\II^{\nn}\right) = \sum\limits_{\substack{P\in V\left(\II^{\mathbf{h}}\colon_R \JJ^{\mathbf{h}}\right),\\ \dim R/P = t}} e(R/P)\cdot l_{R_P}\left({\JJ^{\nn}}_P/{\II^{\nn}}_P\right).
\end{equation*}
For every prime ideal $P$ which appears in the above sum, the $R_P$-module ${\JJ^{\nn}}_P/{\II^{\nn}}_P$ has finite length, and by Theorem \ref{rarepolynomial} it follows that $l_{R_P}\left({\JJ^{\nn}}_P/{\II^{\nn}}_P\right)$ eventually agrees with a numerical polynomial $F_P(\nn)$ of total degree at most $\dim R_P$ and all monomials of highest degree in this polynomial have nonnegative coefficients. Then $e\left(\mathbf{J}^{\nn}/\II^{\nn}\right)$ is eventually given by a numerical polynomial $F(\nn)$ of total degree at most $$\max\left\{\dim R_P \mid P\in V\left(\sqrt{\left(\II^{\mathbf{h}}\colon_R \JJ^{\mathbf{h}}\right)}\right)\;\text{and}\; \dim R/P = t\right\} \leq d-t.$$ As each $e(R/P)>0$ so every monomial of highest degree occuring in $F(\nn)$ has nonnegative coefficient. From \cite[Proposition $2.4$]{cat} or \cite[Proposition $2.1$]{J} we know that there exists a numerical polynomial $G(n)\in \mathbb{Q}[n]$ of degree at most $d-t$ such that $$e\left(\left(J_1\cdots J_r\right)^n/\left(I_1\cdots I_r\right)^n\right) = G(n)$$ for all $n\gg 0$. The last claim then follows from the observation that $F(n,\ldots,n) = G(n)$ for all $n\gg 0$.
\end{proof}

\section{Saturated multi-Rees algebras of ideals with low analytic spread}\label{prooffirst}

The following result subsumes Theorem \ref{mainthm1}. This can be considered as a multi-ideal version of \cite[Theorem $2.1$]{DC4} due to Cutkosky, Herzog and Srinivasan. We recall the definitions of saturated multi-Rees algebra and limit dimension as described in the Introduction and Proposition \ref{constant} respectively.

\begin{theorem}\label{finitely}
Suppose that $(R,m_R)$ is an excellent local domain of dimension $d$ and $I_1,\ldots,I_r,J$ are ideals in $R$. Assume that the analytic spread $\ell({I_1}_P\cdots {I_r}_P) < \dim R_P$ for all $P\in V(J)$. Then the following statements are true.
\begin{enumerate}
\item[$(a)$] The saturated multi-Rees algebra $S_J(I_1,\ldots,I_r)$ of $I_1,\ldots,I_r$ with respect to $J$, is an $\mathbb{N}^r$-graded Noetherian $R$-algebra and also a finitely generated $\NN^r$-graded module over the multi-Rees algebra $\mathcal{R}(I_1,\ldots,I_r)$.
\item[$(b)$] In view of Proposition \ref{constant}, set $$t:= \lim\limits_{\nn\to\infty}\dim_R H^0_J\left(R/\II^{\nn}\right) = \lim\limits_{n\to\infty}\dim_R H^0_J\left(R/\left(I_1\cdots I_r\right)^{n}\right).$$ Then the multigraded multiplicity function $e\left(H^0_J\left(R/\II^{\nn}\right)\right)$ eventually agrees with a numerical polynomial in $\mathbb{Q}[n_1,\ldots,n_r]$  and all monomials of highest degree in this polynomial have nonnegative coefficients. Moreover,
\begin{align*}
\mathrm{total\; degree}\; e\left(H^0_J\left(R/\II^{\nn}\right)\right) &= \mathrm{deg}\; e\left(H^0_J\left(R/\left(I_1\cdots I_r\right)^{n}\right)\right)\\
&\leq \min\left\{d-t-2,\ell(I_1\cdots I_r)+t-1\right\}.
\end{align*}
\end{enumerate}
\end{theorem}

\begin{proof}
$(a)$ Since $(0)\colon_R J^{\infty} = (0)$, we may assume without any loss of generality that $I_1,\ldots,I_r$ are all nonzero ideals. It suffices to show that there are $\NN^r$-graded inclusions $$\mathcal{R}(I_1,\ldots,I_r) \subset S_J(I_1,\ldots,I_r) \subset \overline{\mathcal{R}(I_1,\ldots,I_r)} = \bigoplus\limits_{\nn\in\mathbb{N}^r} \overline{\II^{\nn}\overline{R}}\mathbf{t}^{\nn}$$ and the last equality is a consequence of Lemma \ref{reesint}. Since $$\II^{\nn} \colon_R J^{\infty} \subset \overline{\II^{\nn} \;\overline{R}} \colon_{\overline{R}} J^{\infty}\overline{R}$$ for all $\nn\in\NN^r$, therefore it is enough to show that $$\overline{\II^{\nn} \;\overline{R}} \colon_{\overline{R}} J^{\infty}\overline{R} = \overline{\II^{\nn} \;\overline{R}}$$ for all $\nn\in\NN^r$. We now claim that
\begin{equation}\label{grade1}
\mathrm{grade}(J\overline{R},\overline{R}) := \inf\left\{\mathrm{depth}\;\overline{R}_P \mid P\in V(J\overline{R})\right\} \geq 2.
\end{equation}
Let $P\in V(J\overline{R})$ and $Q = P\cap R$. Then $Q \in V(J)$, and our assumptions imply that $\dim R_Q \geq \ell({I_1}_Q\cdots {I_r}_Q) + 1 \geq 2$. As $R/Q \hookrightarrow \overline{R}/P$ is an integral extension, we obtain that $\dim R/Q = \dim \overline{R}/P$. Hence, since $R$ is an excellent domain, we get $$\dim \overline{R}_P = \dim \overline{R} - \dim \overline{R}/P = \dim R - \dim R/Q = \dim R_Q \geq 2.$$ The ring $\overline{R}$ being normal, satisfies Serre's $S_2$ criterion. Altogether we see that $$\mathrm{depth}\;\overline{R}_P \geq \inf\{2,\dim \overline{R}_P\} \geq 2$$ for all $P\in V(J\overline{R})$ and conclude that $\mathrm{grade}(J\overline{R},R)\geq 2$. Let $A = \overline{\R(I_1,\ldots,I_r)}$ and we further claim that
\begin{equation}\label{grade2}
\mathrm{grade}(JA,A) := \inf\left\{\mathrm{depth}\;A_P \mid P\in V(JA)\right\} \geq 2.
\end{equation}
Again, since $A$ is normal, it satisfies Serre's criterion $S_2$, and it remains to be shown that $\dim A_P \geq 2$ for all $P\in V(JA)$. Let $Q = P\cap R$. Then $Q \in V(J)$. Localizing at $Q$ we may assume that $Q=m_R$. Since $R$ is excellent, we have $$\dim A_P = \dim A - \dim A/P.$$ The excellent property ensures that $A$ is a finite $\R(I_1,\ldots,I_r)$-module and from the discussions in \cite{HHRT}, it follows that $$\dim A = \dim \R(I_1,\ldots,I_r) = \dim R + r.$$ Using \cite[Corollary 3.10]{AJ} we see that $$\dim \dfrac{A}{P} \leq \dim \dfrac{A}{QA} = \dim \dfrac{\R(I_1,\ldots,I_r)}{m_R \R(I_1,\ldots,I_r)} = \ell(I_1\cdots I_r)+r-1.$$ The hypothesis $\ell(I_1\cdots I_r) \leq d -1$ implies that $$\dim A_P \geq (d + r) - (\ell(I_1\cdots I_r) + r -1) \geq 2.$$ From \eqref{grade1} we see that $\mathrm{grade}(J\overline{R},\overline{R})\geq 2$. This implies that $H^0_{J\overline{R}}(\overline{R}) = H^1_{J\overline{R}}(\overline{R}) = 0$ (see \cite[Theorem 6.2.7]{BS}), and hence $$\thickfrac{\overline{\II^{\nn} \;\overline{R}} \colon_{\overline{R}} J^{\infty}\overline{R}}{\overline{\II^{\nn} \;\overline{R}}} \cong H^0_{J\overline{R}}\left(\thickfrac{\overline{R}}{\overline{\II^{\nn} \;\overline{R}}}\right) \cong H^1_{J\overline{R}}\left(\overline{\II^{\nn} \;\overline{R}}\right)$$ for all $\nn\in\NN^r$. Again from \eqref{grade2} we see that $\mathrm{grade}(JA,A)\geq 2$ and using \cite[Theorem 6.2.7]{BS}, we obtain that $$H^1_{JA}(A) = \bigoplus\limits_{\nn\in\NN^r} H^1_{J\overline{R}}\left(\overline{\II^{\nn} \;\overline{R}}\right) = 0.$$ Combining the above results, we conclude that $$\overline{\II^{\nn} \;\overline{R}} \colon_{\overline{R}} J^{\infty}\overline{R} = \overline{\II^{\nn} \;\overline{R}}$$ for all $\nn\in\NN^r$. This proves part $(a)$ of our theorem.

\vspace{0.2cm}
\noindent
$(b)$ Again we may assume without any loss of generality that $I_1,\ldots,I_r$ are nonzero ideals in $R$. It follows from part $(a)$ of this theorem that $$\dfrac{S_J(I_1,\ldots,I_r)}{\mathcal{R}(I_1,\ldots,I_r)} = \bigoplus\limits_{\nn\in\NN^r}\dfrac{\II^{\nn} \colon_R J^{\infty}}{\II^{\nn}} = \bigoplus\limits_{\nn\in\NN^r} H^0_J\left(\dfrac{R}{\II^{\nn}}\right)$$ is a finitely generated $\NN^r$-graded module over $\mathcal{R}(I_1,\ldots,I_r)$. We know from Proposition \ref{constant} that there exists a vector $\mathbf{h}\in\mathbb{N}^r$ such that $$\mathrm{Ann}_R \left(H^0_J\left(R/\II^{\mathbf{n}}\right)\right) = \mathrm{Ann}_R \left(H^0_J\left(R/\II^{\mathbf{h}}\right)\right)=: Q$$ for all $\mathbf{n}\geq \mathbf{h}$. Then $t = \dim R/Q$, $\sqrt{I+J} \subseteq \sqrt{Q}$ and there exist equalities $H^0_J\left(R/\II^{\mathbf{n}}\right) = H^0_{\sqrt{Q}}\left(R/\II^{\mathbf{n}}\right)$ for all $\mathbf{n}\geq \mathbf{h}$. From the associativity formula we get
\begin{align}\label{mainsum}
e\left(H^0_J\left(R/\II^{\nn}\right)\right) &= \sum\limits_{P\in V(Q),\;\dim R/P = t}e(R/P)\cdot l_{R_P}\left(H^0_{P_P}\left(R_P/{\II^{\nn}}_P\right)\right)
\end{align}
for all $\nn\geq \mathbf{h}$. For every prime ideal $P$ appearing in the above sum \eqref{mainsum}, we have that $\bigoplus_{\nn\in\NN^r}H^0_{P_P}\left(R_P/{\II^{\nn}}_P\right)$ is a finitely generated $\NN^r$-graded module over the multi-Rees algebra $\mathcal{R}({I_1}_P,\ldots,{I_r}_P)$ and there exists an integer $c>0$ such that $P_P^c$ annihilates the module $\bigoplus_{\nn\in\NN^r}H^0_{P_P}\left(R_P/{\II^{\nn}}_P\right)$. So it can be regarded as a finitely generated $\NN^r$-graded module over $\mathcal{R}({I_1}_P,\ldots,{I_r}_P)/P_P^c\mathcal{R}({I_1}_P,\ldots,{I_r}_P)$ and using \cite[Corollary 3.10]{AJ} we notice that
\begin{equation}\label{atmostdim}
\dim \dfrac{\mathcal{R}({I_1}_P,\ldots,{I_r}_P)}{P_P^c \mathcal{R}({I_1}_P,\ldots,{I_r}_P)} -r = \ell\left({I_1}_P,\ldots,{I_r}_P\right) - r = \ell\left({I_1}_P\cdots {I_r}_P\right)-1.
\end{equation}
From \cite[Lemma $1.1$ and Theorem $4.1$]{HHRT} and statement \eqref{atmostdim} we know that there exists a numerical polynomial $F_P(\nn) \in\mathbb{Q}[n_1,\ldots,n_r]$ of total degree at most $\ell\left({I_1}_P\cdots {I_r}_P\right)-1$ such that $$l_{R_P}\left(H^0_{P_P}\left(R_P/{\II^{\nn}}_P\right)\right)=F_P(\nn)$$ for all $\nn\gg\mathbf{0}$ and all monomials of highest degree in this polynomial have nonnegative coefficients. We finally conclude from \eqref{mainsum} that $e\left(H^0_J\left(R/\II^{\nn}\right)\right)$ eventually agrees with a numerical polynomial $F(\nn):= \sum\limits_{P\in V(Q),\;\dim R/P=t} F_P(\nn)$, of total degree at most
\begin{align*}
&\max\left\{\ell\left({I_1}_P\cdots {I_r}_P\right)-1 \mid P\in V(Q),\;\dim R/P = t\right\}\\
&\leq \max\left\{\dim R_P - 2 \mid P\in V(Q),\;\dim R/P = t\right\}\\
&\leq d-t-2
\end{align*}
and all monomials of highest degree in this polynomial have nonnegative coefficients. From the proof of \cite[Theorem $2.1(b)$]{DC4} we know that there exists a numerical polynomial $G(n)\in \mathbb{Q}[n]$ of degree at most $\ell(I_1\cdots I_r) + t -1$ such that $$e\left(H^0_J\left(R/\left(I_1\cdots I_r\right)^{n}\right)\right) = G(n)$$ for all $n\gg 0$. The last claim is a consequence of the observation that $F(n,\ldots,n) = G(n)$ for all $n\gg 0$.
\end{proof}

A refinement of Theorem \ref{finitely} is obtained when $J=m_R$ and this observation is recorded below. It generalises and improves statement $(0.6)$ of \cite{DC4}.

\begin{corollary}\label{smallpolynomial}
Suppose that $(R,m_R)$ is an excellent local domain of dimension $d$ and $I_1,\ldots,I_r$ are ideals in $R$. Assume that the analytic spread $\ell({I_1}\cdots {I_r}) < d$. Then the saturated multi-Rees algebra $S_{m_R}(I_1,\ldots,I_r)$ is an $\mathbb{N}^r$-graded Noetherian $R$-algebra and also a finitely generated $\NN^r$-graded module over the multi-Rees algebra $\mathcal{R}(I_1,\ldots,I_r)$. The multigraded length function $l_R\left(H^0_{m_R}\left(R/\II^{\nn}\right)\right)$ eventually agrees with a numerical polynomial in $\mathbb{Q}[n_1,\ldots,n_r]$ and all monomials of highest degree in this polynomial have nonnegative coefficients. In particular,
$$\mathrm{total\; degree}\; l_R\left(H^0_{m_R}\left(R/\II^{\nn}\right)\right) = \mathrm{deg}\; l_R\left(H^0_{m_R}\left(R/\left(I_1\cdots I_r\right)^{n}\right)\right) \leq \ell(I_1\cdots I_r)-1.$$
\end{corollary}

We collect some known facts about analytic spread of ideals.

\begin{remarks}
Suppose that $(R,m_R)$ is a Noetherian local domain of dimension $d$ and $I_1,\ldots,I_r$ are nonzero ideals in $R$.
\begin{enumerate}
\item[$(i)$] If $\ell\left(I_1\cdots I_r\right)<d$ then it follows from \cite[Corollary $3.12$]{AJ} that each $\ell(I_j)<d$. However the converse need not be true. Consider the ring $R=K[X,Y,Z]$ with ideals $I_1 = (X,Y)$ and $I_2 = (Y,Z)$. In this case, $\ell(I_1) = 2$, $\ell(I_2)=2$ whereas $\ell(I_1I_2)=3$.
\item[$(ii)$] If each $\ell(I_j)\leq (d-2)/r+1$ then it again follows from \cite[Corollary $3.12$]{AJ} that $\ell(I_1\cdots I_r)\leq \sum_{j=1}^r \ell(I_j) - r +1 \leq d-1$.
\item[$(iii)$] If at least one of the ideals $I_j$ is principal then by possibly permuting the indices we may assume that $I_1,\ldots,I_s$ are nonzero principal ideals for some $s\leq r$. Then \cite[Corollary $3.8$]{AJ} shows that $\ell(I_1\cdots I_r) = \ell(I_{s+1}\cdots I_r)$.
\end{enumerate}
\end{remarks}

\section{Saturated multi-Rees algebras of monomial ideals}

Throughout this section we shall assume that $R=K[X_1,\ldots,X_d]$ is a standard graded polynomial ring in $d$ variables over a field $K$ and $m_R=(X_1,\ldots,X_d)$ is the graded maximal ideal of $R$, unless otherwise specified. For any $\mathbf{a}= (a_1,\ldots,a_d)\in\NN^d$ we shall denote the monomial $X_1^{a_1}\cdots X_d^{a_d}$ by $\mathbf{X}^{\mathbf{a}}$.

\begin{lemma}\label{monomialideal}
 Suppose that $I$ is a monomial ideal and $J$ is any ideal in $R$. Then $I\colon_R J^{\infty}$ is also a monomial ideal in $R$.
\end{lemma}
\begin{proof}
 Since $I$ is a monomial ideal, there exists a primary decomposition $I = \bigcap_{i=1}^m Q_i$, where each $Q_i$ is a monomial ideal generated by pure powers of variables \cite[Theorem 1.3.1]{HH}. By \cite[Proposition $3.13(a)$]{eisenbud} we have $$I\colon_R J^{\infty} = \bigcap\limits_{\substack{i=1,\\J\not\subset \sqrt{Q_i}}}^m Q_i.$$ Therefore $I\colon_R J^{\infty}$ is a monomial ideal as finite intersections of monomial ideal is again a monomial ideal \cite[Proposition 1.2.1]{HH}.
\end{proof}

\begin{lemma}\label{pigeonhole}
Suppose that $S$ is a commutative ring with identity and $I,J_1,J_2$ are ideals in $S$. Then $$I \colon_S \left(J_1+J_2\right)^{\infty} = \left(I \colon_S J_1^{\infty}\right)\cap \left(I \colon_S J_2^{\infty}\right).$$
\end{lemma}
\begin{proof}
This proof is fairly simple and only involves elementary arguments.
\end{proof}
For any nonzero polynomial $f\in K[X_1,\ldots,X_d]$ there is a unique writing $f = \sum_{\nn\in\NN^d}a_{\nn}\X^{\nn}$ where the coefficients $a_{\nn}\in K$ and $\X^{\nn}:=X_1^{n_1}\cdots X_d^{n_d}$. Define the monomial ideal $\mathcal{S}_f := \left(\X^{\nn} \mid a_{\nn}\neq 0 \right)$. In other words, $\mathcal{S}_f$ is the ideal generated by the monomials occurring in the expression of $f$ with nonzero coefficients. Let $>_{\mathrm{mon}}$ be any \emph{monomial ordering} \cite[Definition 1]{cox} on $R$. The \emph{leading term} \cite[Definition 7]{cox} of $f$ is $\mathrm{LT}(f):= a_{\mm}\X^{\mm}$ where $\mm = \max\{\nn\in\NN^d \mid a_{\nn}\neq 0\}$ and the maximum is taken with respect to $>_{\mathrm{mon}}$.

\begin{lemma}\label{symbol2}
Let $f$ be a nonzero polynomial in $R$ and let $\mathcal{S}_{f}$ be the ideal generated by the monomials occurring in the expression of $f$ with nonzero coefficients. Then $$I \colon_R (f)^{\infty} = I \colon_R \mathcal{S}_{f}^{\infty}$$ for all monomial ideals $I$ in $R$.
\end{lemma}
\begin{proof}
We may write $f = \sum_{i=1}^{t}c_i\X^{\mathbf{v}_i}$, where $\X^{\mathbf{v}}:= X_1^{v_{i_1}}\cdots X_d^{v_{i_d}}$ are the distinct monomials occurring in the expression of $f$ with nonzero coefficients $c_i\in K$. Fix a monomial ideal $I$ in $R$ and we shall prove that $I \colon_R (f)^{\infty} = I \colon_R \mathcal{S}_{f}^{\infty}$ by induction on $t$. The base case $t=1$ is obvious since $(f)=\mathcal{S}_f$. Hence assume that $t\geq 2$. Let $>_{\mathrm{mon}}$ be any monomial ordering on $R$. We first claim that 
\begin{equation}\label{leadingterm}
 I \colon_R (f)^{\infty} \subseteq I \colon_R \left(\mathrm{LT}(f)\right)^{\infty}.
\end{equation}
We know from Lemma \ref{monomialideal} that $I \colon_R (f)^{\infty}$ is a monomial ideal. Suppose that $\X^{\mathbf{a}}:=X_1^{a_1}\cdots X_d^{a_d} \in I \colon_R (f)^{\infty}$. Then there exists an integer $r>0$ such that $\X^{\mathbf{a}}f^r \in I$ and therefore $\mathrm{LT}\left(\X^{\mathbf{a}}f^r\right) = \X^{\mathbf{a}}\left(\mathrm{LT}(f)\right)^r \in I$. Thus $\X^{\mathbf{a}} \in I \colon_R \left(\mathrm{LT}(f)\right)^{\infty}$ and proves our claim \eqref{leadingterm}. Now observe that
\begin{align*}
 I \colon_R (f)^{\infty} &= \left(I \colon_R (f)^{\infty}\right)\cap \left(I \colon_R \left(\mathrm{LT}(f)\right)^{\infty}\right) && (\text{claim}\;\eqref{leadingterm})\\
 &= \left(I \colon_R \left(f-\mathrm{LT}(f)\right)^{\infty}\right)\cap \left(I \colon_R \left(\mathrm{LT}(f)\right)^{\infty}\right) && (\text{Lemma}\;\ref{pigeonhole})\\
 &= \left(I \colon_R \mathcal{S}_{\left(f-\mathrm{LT}(f)\right)}^{\infty}\right)\cap \left(I \colon_R \left(\mathrm{LT}(f)\right)^{\infty}\right) && (\text{induction hypothesis})\\
 &= I \colon_R \mathcal{S}_f^{\infty} && (\text{Lemma}\;\ref{pigeonhole}).
\end{align*}
\end{proof}

For every subset $F\subseteq \{1,\ldots,d\}$, define the map $\pi_F \colon R \to R$ by $\pi_F(X_i) = 1$ if $i\in F$, and $\pi_F(X_i) = X_i$ if $i\not\in F$. Note that $\pi_F$ is a ring homomorphism.

\begin{lemma}\label{symbol3}
 Suppose that $\X^{\mathbf{c}}:=X_1^{c_1}\cdots X_d^{c_d}$ is a monomial in $R$ and consider the subset $F = \{i \mid c_i>0\}\subseteq\{1,\ldots,d\}$. Then $$I \colon_R \left(\X^{\mathbf{c}}\right)^{\infty} = \pi_F(I)$$ for all monomial ideals $I$ in $R$.
\end{lemma}
\begin{proof}
 Fix a monomial ideal $I$ in $R$. Then both $I \colon_R \left(\X^{\mathbf{c}}\right)^{\infty}$ and $\pi_F(I)$ are monomial ideals. By possibly renaming the variables we may assume that $c_1>0,\ldots, c_s>0$ and $c_{s+1}=\cdots = c_d = 0$, i.e., $F=\{1,\ldots,s\}$ for some integer $s\leq d$. Consider the two multiplicatively closed subsets $$W:= \left\{X_1^{c_1n}\cdots X_s^{c_sn}\mid n\geq 0\right\} \;\subseteq\; \left\{X_1^{n_1}\cdots X_s^{n_s}\mid n_1,\ldots,n_s\geq 0\right\}=:W^{\prime}$$ in $R$. Note that $W^{-1}R = {W^{\prime}}^{-1}R$ and we have
 \begin{align*}
 I \colon_R \left(X_1^{c_1}\cdots X_s^{c_s}\right)^{\infty} &= I\left(W^{-1}R\right)\cap R\\
 &= I\left({W^{\prime}}^{-1}R\right)\cap R\\
 &= \pi_F(I)\left({W^{\prime}}^{-1}R\right)\cap R\\
 &= \pi_F(I) \colon_R \left(X_1^{c_1}\cdots X_s^{c_s}\right)^{\infty}\\
 &= \pi_F(I).
 \end{align*}
\end{proof}

\begin{example}
We pictorially demonstrate how to compute the saturation of the monomial ideal $I = (x^3y,x^2y^2,xy^3)$ in the polynomial ring $R=k[x,y]\mathrm{:}$

\begin{center}
\begin{tikzpicture}[scale=0.35]
\draw[->] (0,0)--(6.5,0) node[right]{$x$};
\draw[->] (0,0)--(0,6.5) node[above]{$y$};
\draw[,blue] (1,6.5)--(1,3)--(2,3)--(2,2)--(3,2)--(3,1)--(6.5,1);
    \draw[fill=lightCyan, draw=none] (1,6.5)--(1,3)--(2,3)--(2,2)--(3,2)--(3,1)--(6.5,1)--(6.5,6.5)--(1,6.5);
   \draw[thick,blue] (1,6.5)--(1,3)--(2,3)--(2,2)--(3,2)--(3,1)--(6.5,1);
 \foreach \x in {0,...,6}
 \foreach \y in {0,...,6}
 \fill[black,fill=black] (\x,\y) circle (2pt);
\node[blue] at (3,1) {\small{$\bullet$}};
\node[blue] at (2,2) {\small{$\bullet$}};
\node[blue] at (1,3) {\small{$\bullet$}};
\end{tikzpicture}
\begin{tikzpicture}[scale=0.35]
\draw[->] (0,0)--(6.5,0) node[right]{$x$};
\draw[->] (0,0)--(0,6.5) node[above]{$y$};
\draw[fill=lightCyan, draw=none] (6.5,1)--(6.5,6.5)--(0,6.5)--(0,1);
\draw[fill=deepCyan, draw=none] (1,6.5)--(1,3)--(2,3)--(2,2)--(3,2)--(3,1)--(6.5,1)--(6.5,6.5)--(1,6.5);
\draw[thick,blue] (6.5,1)--(0,1)--(0,6.5);
\foreach \x in {0,...,6}
\foreach \y in {0,...,6}
\fill[black,fill=black] (\x,\y) circle (2pt);
\node[blue] at (0,1) {\small{$\bullet$}};
\end{tikzpicture}
\begin{tikzpicture}[scale=0.35]
\draw[->] (0,0)--(6.5,0) node[right]{$x$};
\draw[->] (0,0)--(0,6.5) node[above]{$y$};
\draw[fill=lightCyan, draw=none] (1,6.5)--(1,0)--(6.5,0)--(6.5,6.5)--(1,6.5);
\draw[fill=deepCyan, draw=none] (1,6.5)--(1,3)--(2,3)--(2,2)--(3,2)--(3,1)--(6.5,1)--(6.5,6.5)--(1,6.5);
\draw[thick,blue] (1,6.5)--(1,0)--(6.5,0);
\foreach \x in {0,...,6}
\foreach \y in {0,...,6}
\fill[black,fill=black] (\x,\y) circle (2pt);
\node[blue] at (1,0) {\small{$\bullet$}};
\end{tikzpicture}
\begin{tikzpicture}[scale=0.35]
\draw[->] (0,0)--(6.5,0) node[right]{$x$};
\draw[->] (0,0)--(0,6.5) node[above]{$y$};
\draw[fill=lightCyan, draw=none] (6.5,1)--(6.5,6.5)--(1,6.5)--(1,1);
\draw[fill=deepCyan, draw=none] (1,6.5)--(1,3)--(2,3)--(2,2)--(3,2)--(3,1)--(6.5,1)--(6.5,6.5)--(1,6.5);
\draw[thick,blue] (1,6.5)--(1,1)--(6.5,1);
\foreach \x in {0,...,6}
\foreach \y in {0,...,6}
\fill[black,fill=black] (\x,\y) circle (2pt);
\node[blue] at (1,1) {\small{$\bullet$}};
\end{tikzpicture}
\end{center}
\end{example}

The following result is a straightforward generalisation of \cite[Corollary 1.3]{HHT}. Similar results are also present in \cite[Chapter IV]{fields}.

\begin{proposition}\label{semigroup}
Suppose that ${I}_{ij}$ is a monomial ideal in $R$ for each $i=1,\ldots,r$ and for each $j=1,\ldots,s$. Then the intersection algebra $$\bigoplus\limits_{(n_1,\ldots,n_r)\in\NN^r}\left(\bigcap\limits_{j=1}^s {I}_{1j}^{n_1}\cdots {I}_{rj}^{n_r} \right)t_1^{n_1}\cdots t_r^{n_r}$$ is an $\mathbb{N}^r$-graded Noetherian $R$-algebra.
\end{proposition}
\begin{proof}
For every $j\in\{1,\ldots,s\}$, let $A_j := \bigoplus\limits_{(n_1,\ldots,n_r)\in\NN^r}{I}_{1j}^{n_1}\cdots {I}_{rj}^{n_r}t_1^{n_1}\cdots t_r^{n_r}$ be the multi-Rees algebra of the monomial ideals ${I}_{1j},\ldots, {I}_{rj}$ in $R$ and let $A := \bigcap_{j=1}^s A_j$ where the intersection is taken inside $R[t_1,\ldots,t_r]$. Then these algebras can be viewed as semigroup rings, i.e., $A = K[H]$ and $A_j = K[H_j]$ for all $j=1,\ldots,s$ and $H,H_1,\ldots,H_s$ are affine semigroups in $\ZZ^{d+r}$. The affine semigroups $H_j$ are finitely generated since the multi-Rees algebras $A_j$ are finitely generated and $H = \bigcap_{j=1}^s H_j$ since $A=\bigcap_{j=1}^s A_j$. \cite[Corollary 1.2]{HHT} implies that $H$ is finitely generated and consequently $A$ is finitely generated as an $R$-algebra.
\end{proof}

For the next result we need to introduce some notations which we borrow from \cite{KW}. A map $\sigma \colon \NN^r \to \mathbb{Q}$ is said to be \emph{periodic} if there exists an $\alpha\in\NN$ such that $$\sigma(\nn + \alpha\cdot \mathbf{e}_i) = \sigma(\nn)$$ for all $\nn\in\NN^r$ and $i=1,\ldots,r$. A \emph{quasi-polynomial of degree $d$} (over $\mathbb{Q}$) is a function $f \colon \NN^r \to \mathbb{Q}$ such that for all $\nn\in\NN^r$ there exists an expression $$f(\nn) = \sum\limits_{\mathbf{i}\in\NN^r,\;|\mathbf{i}|\leq d}\sigma_{\mathbf{i}}(\nn)\nn^{\mathbf{i}}$$ where $\sigma_{\mathbf{i}} \colon \NN^r \to \mathbb{Q}$ are periodic functions for all $\mathbf{i}\in\NN^r$ with $|\mathbf{i}|\leq d$. A \emph{piecewise quasi-polynomial} (over $\mathbb{Q}$) is a function $f \colon \NN^r \to \mathbb{Q}$ such that there exists a finite partition of $\NN^r = \bigcup_{i=1}^t \left(P_i \cap \NN^r\right)$ with $P_i$ polyhedra (which may not all be full-dimensional) and there exist quasi-polynomials $f_i$ such that $$f(\nn) = f_i(\nn)$$ for all $\nn\in P_i\cap \NN^r$ and $i=1,\ldots,t$.

A family of ideals $\mathcal{I}=\{I_{\nn}\}_{\nn\in\NN^r}$ in $R$ is said to be a \emph{Noetherian $\NN^r$-graded family} if $I_{\mathbf{0}}=R$ and $\bigoplus_{\nn\in\NN^r}I_{\nn}$ is a Noetherian $\NN^r$-graded $R$-algebra.

\begin{lemma}\cite[Proposition $1$]{das2}\label{piecewise}
Suppose that $\mathcal{F}=\{F_{\mathbf{n}}\}_{\mathbf{n}\in\mathbb{N}^r}$ and $\mathcal{G}=\{G_{\mathbf{n}}\}_{\mathbf{n}\in\mathbb{N}^r}$ be two Noetherian $\mathbb{N}^r$-graded families of $R$ by monomial ideals. For every $\mathbf{n}\in\mathbb{N}^r$, consider the set $$S_{\mathbf{n}} = \left\{\mathbf{a}\in\mathbb{N}^d \mid \mathbf{X}^{\mathbf{a}} \in F_{\mathbf{n}}\;\text{and}\;\mathbf{X}^{\mathbf{a}} \not\in G_{\mathbf{n}}\right\}.$$ Assume that $\#S_{\mathbf{n}}$ is finite for all $\mathbf{n}\gg 0$. Then there exists a piecewise quasi-polynomial $F \colon \mathbb{N}^r \to \mathbb{N}$ such that $F(\mathbf{n}) = \#S_{\mathbf{n}}$ for all $\mathbf{n}\gg \mathbf{0}$.
\end{lemma}

\subsection{Main result}\label{proofsecond}

\begin{proof}[Proof of Theorem \ref{mainthm2}]
We will only consider the case when $L\neq (0)$. It is well-known that every irreducible monomial ideal is of the form $\left(X_{i_1}^{a_1},\ldots,X_{i_s}^{a_s}\right)$ for some integers $s\leq d$, $i_1<\cdots<i_s$ and $a_i\geq 1$. By possibly renaming the variables we may assume that $L = \left(X_1^{a_1},\ldots,X_s^{a_s}\right)$. Consider the natural map $\varphi\colon R \to R/L=:S$ and write $I_1^{\prime} = \varphi^{-1}(I_1),\ldots,I_r^{\prime} = \varphi^{-1}(I_r), J^{\prime} = \varphi^{-1}(J)$. Note that $I_1^{\prime},\ldots,I_r^{\prime}$ are monomial ideals in $R$. It follows from combining Lemma \ref{pigeonhole}, Lemma \ref{symbol2} and Lemma \ref{symbol3} that there exist subsets $\mathcal{A}_1,\ldots,\mathcal{A}_l \subseteq \{1,\ldots,d\}$, depending only on $J^{\prime}$, such that $I^{\prime}\colon_R {J^{\prime}}^{\infty} = \bigcap_{i=1}^l \pi_{\mathcal{A}_i}\left(I^{\prime}\right)$ for all monomial ideals $I^{\prime}$ in $R$. Therefore
\begin{align*}
{I_1}^{n_1}\cdots {I_r}^{n_r} \colon_S J^{\infty} &= \left({I_1^{\prime}}^{n_1}\cdots {I_r^{\prime}}^{n_r} + L\right) \colon_R {J^{\prime}}^{\infty} \mod L\\
&= \bigcap\limits_{1\leq i\leq l} \left(\pi_{\mathcal{A}_i}\left(I_1^{\prime}\right)^{n_1}\cdots\pi_{\mathcal{A}_i}\left(I_r^{\prime}\right)^{n_r} + \pi_{\mathcal{A}_i}\left(L\right)\right) \mod L\\
&= \left(\bigcap\limits_{\substack{1\leq i\leq l,\\ \mathcal{A}_i\subseteq \{s+1,\ldots,d\}}} \pi_{\mathcal{A}_i}\left(I_1^{\prime}\right)^{n_1}\cdots\pi_{\mathcal{A}_i}\left(I_r^{\prime}\right)^{n_r}\right) + L \mod L
\end{align*}
for all non-negative integers $n_1,\ldots,n_r$. In the last equality above, we have used that $L$ is an irreducible monomial ideal. We know from Proposition \ref{semigroup} that $$\bigcap\limits_{\substack{1\leq i\leq l,\\ \mathcal{A}_i\subseteq \{s+1,\ldots,d\}}} \pi_{\mathcal{A}_i}\left(I_1^{\prime}\right)^{n_1}\cdots\pi_{\mathcal{A}_i}\left(I_r^{\prime}\right)^{n_r}$$ is a Noetherian $\NN^r$-graded family of monomial ideals in $R$ and this proves the first assertion of our theorem.

When $J=m_S$ then $J^{\prime}=m_R$, the homogeneous maximal ideal of $R$, and $\mathcal{A}_i = \{i\}$ for all $i=1,\ldots,d$. Moreover, the length function $l_S\left(H^0_{m_S}\left(S/I_1^{n_1}\cdots I_r^{n_r}\right)\right)$ equals $$\#\left\{(a_1,\ldots,a_d)\in\mathbb{N}^d \mid X_1^{a_1}\cdots X_d^{a_d}\in \left(\bigcap\limits_{i=s+1}^d \pi_{\{i\}}\left(I_1^{\prime}\right)^{n_1}\cdots\pi_{\{i\}}\left(I_r^{\prime}\right)^{n_r}\right)\setminus {I_1^{\prime}}^{n_1}\cdots {I_r^{\prime}}^{n_r}\right\}.$$ The second assertion now follows directly from Lemma \ref{piecewise} and Theorem \ref{rightbound}.
\end{proof}

The following example from \cite[Example $1$]{das2} illustrates the piecewise polynomial type behaviour of the multigraded length function $l_R\left(H^0_{m_R}\left(R/I_1^{n_1}\cdots I_r^{n_r}\right)\right)$.

\begin{example}
Suppose that $R=K[X,Y,Z]$ with homogeneous maximal ideal $m_R=(X,Y,Z)$. Consider the monomial ideals $I_1 = (X,Y)$ and $I_2 = (Y,Z)$. Then we have a formula $$l_R\left(H^0_{m_R}\left(R/I_1^{n_1}I_2^{n_2}\right)\right) = \begin{cases}
\frac{n_1(n_1+1)(3n_2-n_1+1)}{6} &n_2\geq n_1,\\
\frac{n_1(n_2+1)(3n_1-n_2+1)}{6} &n_1\geq n_2
\end{cases}$$
for all $n_1,n_2\in \NN$.
\end{example}

The following corollary gives the combined statements of \cite[Theorem $3.2$]{HHT} and \cite[Corollary $2.5$]{J} when $L=(0)$ and $J$ is a monomial ideal.

\begin{corollary}\label{maincor}
Consider the setup in Theorem \ref{mainthm2}, i.e., $R = K[X_1,\ldots,X_d]$ and $S=R/L$ where $L\subseteq R$ is either the zero ideal or an irreducible monomial ideal. Suppose that $I$ is an ideal generated by monomials in $S$ and $J$ is any ideal in $S$. Then the saturated Rees algebra $S_J(I)$ of $I$ with respect to $J$ is Noetherian. Moreover, there exist an integer $g>0$ and numerical polynomials $P_0,\ldots,P_{g-1} \in \mathbb{Q}[n]$ such that the following conditions hold:
\begin{enumerate}
 \item[$(i)$] $e\left(H^0_J\left(S/I^{ng+r}\right)\right) = P_{r}(n)$ for all $n\gg 0$ and $r=0,\ldots, g-1$. In other words, the multiplicity function $e\left(H^0_J\left(S/I^{n}\right)\right)$ is of quasi-polynomial type of period $g$.
 \item[$(ii)$] The limit dimension $t:= \lim\limits_{n\to\infty}\dim_S H^0_J\left(S/I^n\right)$ exists.
 \item[$(iii)$] $\deg P_r = \deg P_0 \leq \dim S - t$ for all $r=1,\ldots,g-1$ and all of them have the same leading coefficient. In particular, the limit $$\lim\limits_{n\to\infty}\dfrac{e\left(H^0_J\left(S/I^{n}\right)\right)}{n^{\dim S -t}}$$ exists and is a rational number.
\end{enumerate}
\end{corollary}
\begin{proof}
Consider the graded inclusion of graded $R$-algebras $$\mathcal{R}(I) = \bigoplus_{n\in\NN}I^n \subseteq \bigoplus_{n\in\NN}\left(I^n\colon_S J^{\infty}\right) = S_J(I).$$ Theorem \ref{mainthm2} says that $S_J(I)$ is a finitely generated $S$-algebra. If $\mathrm{grade}(I,S)=0$ then $I$ is a nilpotent ideal and we are done. Hence assume that $\mathrm{grade}(I,S)>0$ and the remaining assertions now follow from \cite[Theorem $3.2$]{J}.
\end{proof}

Corollary \ref{maincor} extends the class of ideals for which the epsilon multiplicity exists as a limit and is rational. The following example suggests that the assumptions on $L$ in Corollary \ref{maincor} are necessary to guarantee finite generation of the saturated Rees algebra $S_J(I)$.

\begin{example}\cite[Example $2.2$]{DC4}
Suppose that $S = K[X,Y,Z]/(XY,XZ)$, $I = (X,Y)S$ and $m_S = (X,Y,Z)S$. Then the saturated Rees algebra $$S_{m_S}(I) = \bigoplus\limits_{n\in\NN}\left(I^n\colon_S m_S^{\infty}\right) = \bigoplus\limits_{n\in\NN}\left(X,Y^n\right)S$$ is not a finitely generated $S$-algebra.
\end{example}

In view of the previous results, we now pose two natural questions.

\begin{question}
Suppose that $R = R_0[X_1,\ldots,X_d]$ is a standard graded polynomial ring in $d$ variables over an Artinian local ring $(R_0, m_{R_0})$ and $m_R:= m_{R_0} + (X_1,\ldots,X_d)$ is the graded maximal ideal of $R$. Suppose that $I\subseteq R$ is a graded ideal generated by monomial terms of the form $c_{i_1,\ldots,i_r}X_1^{i_1}\cdots X_d^{i_d}$ where $i_1,\ldots,i_d \in \NN$ and $c_{i_1,\ldots,i_r} \in R_0$. Does the epsilon multiplicity $\varepsilon(I)$ exist as a limit and is it rational? More generally, is the saturated Rees algebra $S_{m_R}(I)$ of $I$ Noetherian? Can we relate $\varepsilon(I)$ to the volume of a certain region, in the sense of \cite{JM}?
\end{question}

\begin{question}
Suppose that $(R,m_R)$ is an analytically unramified local ring of dimension $d$ and $I_1,\ldots,I_r$ are ideals in $R$ such that the analytic spread $\ell(I_1\cdots I_r)=d$. Is it possible to define `mixed' versions of epsilon multiplicity, in the sense of \cite{CSS}, by studying the function $$P(n_1,\ldots,n_r) = \lim\limits_{n\to\infty}\dfrac{l_R\left(H^0_{m_R}\left(R/I_1^{nn_1}\cdots I_r^{nn_r}\right)\right)}{n^d} \quad \text{for $n_1,\ldots,n_r\gg 0$}.$$
\end{question}

\section{Saturated multi-Rees algebras in dimension two}\label{dimensiontwo}

\begin{proposition}\label{dimension2}
Suppose that $(R,m_R)$ is a Noetherian local unique factorization domain of dimension two and $I_1,\ldots,I_r$ are ideals in $R$. Then the saturated multi-Rees algebra $S_{m_R}(I_1,\ldots,I_r)$ is a standard $\NN^r$-graded Noetherian $R$-algebra. If the analytic spread $\ell(I_1\cdots I_r)=2$ then there exists a numerical polynomial $F\in\mathbb{Q}[n_1,\ldots,n_r]$ of total degree $2$ such that $$l_R\left(H^0_{m_R}\left(R/\II^{\nn}\right)\right) = F(\nn)$$ for all $\nn\gg\mathbf{0}$ and all monomials of highest degree in $F$ have nonnegative coefficients. Otherwise if the analytic spread $\ell(I_1\cdots I_r)\leq 1$ then this length function is zero.
\end{proposition}
\begin{proof}
Let $j$ be an integer such that $1\leq j\leq r$. Without any loss of generality we may assume that $I_j\neq (0)$ for all such $j$. If $I_j$ is not a principal ideal then we can write $I_j = (f_j)Q_j$, where $Q_j$ is an $m_R$-primary ideal and $f_j$ is the g.c.d of the nonzero elements of $I_j$. Observe that if $I_j$ is $m_R$-primary then $f_j=1$ and $Q_j = I_j$. Otherwise, if $I_j = (f_j)$ is a principal ideal then we set $Q_j = R$. For any $\nn\in\NN^r$ we can verify that $$\II^{\nn}\colon_R m_R^{\infty} = (\mathbf{f}^{\nn})$$ where $\mathbf{f}^{\nn}$ denotes the product $f_1^{n_1}\cdots f_r^{n_r}$. Thus we see that $$S_{m_R}(I_1,\ldots,I_r) = \bigoplus\limits_{\nn\in\NN^r}\left(\II^{\nn}\colon_R m_R^{\infty}\right)\mathbf{t}^{\nn} = \bigoplus\limits_{\nn\in\NN^r}(\mathbf{f}^{\nn}) \mathbf{t}^{\nn}$$ is a standard $\NN^r$-graded Noetherian $R$-algebra. Moreover, there exist isomorphisms $$H^0_{m_R}\left(R/\II^{\nn}\right) = \left(\mathbf{f}^{\nn}\right)/\left(\mathbf{f}^{\nn}\right)\mathbf{Q}^{\nn} \cong R/\mathbf{Q}^{\nn}$$ where $\mathbf{Q}^{\nn}$ denotes the product $Q_1^{n_1}\cdots Q_r^{n_r}$. The remaining assertions now follow.
\end{proof}

\section{Acknowledgements}

The authors would like to thank their Ph.D. advisors, Prof. Steven Dale Cutkosky, Prof. Tony Joseph Puthenpurakal and Prof. Jugal Kishore Verma. Additionally, they thank Prof. Manoj Kummini for suggesting a simpler proof of Lemma \ref{symbol2}. Major part of this work was carried out when both the authors were postdocs at the Chennai Mathematical Institute, India. The authors are grateful to the Infosys Foundation for providing partial financial support during their stay at CMI.

\bibliographystyle{alpha}
\bibliography{main}

\begin{thebibliography}{{Woo}15}

\bibitem[Ama76]{JO}
J.~O. Amao.
\newblock On a certain {H}ilbert polynomial.
\newblock {\em Journal of the London Mathematical Society}, 2(1):13--20, 1976.

\bibitem[BAM20]{AJ}
C.~Bivi{\`a}-Ausina and J.~Monta{\~n}o.
\newblock Analytic spread and integral closure of integrally decomposable
  modules.
\newblock {\em Nagoya Mathematical Journal}, pages 1--26, 2020.

\bibitem[BF00]{fields}
J.~Bruce~Fields.
\newblock {\em Length functions determined by killing powers of several ideals
  in a local ring}.
\newblock PhD thesis, University of Michigan, 2000.

\bibitem[BS12]{BS}
M.~P. Brodmann and R.~Y. Sharp.
\newblock {\em Local cohomology: an algebraic introduction with geometric
  applications}, volume 136.
\newblock Cambridge university press, 2012.

\bibitem[CHS10]{DC4}
S.~D. Cutkosky, J.~Herzog, and H.~Srinivasan.
\newblock Asymptotic growth of algebras associated to powers of ideals.
\newblock {\em Mathematical Proceedings of the Cambridge Philosophical
  Society}, 148(1):55--72, 2010.

\bibitem[CHST05]{DC3}
S.~D. Cutkosky, H.~T. H{\`a}, H.~Srinivasan, and E.~Theodorescu.
\newblock Asymptotic behavior of the length of local cohomology.
\newblock {\em Canadian Journal of Mathematics}, 57(6):1178--1192, 2005.

\bibitem[{Ciu}15]{cat}
C.~{Ciuperc\u{a}}.
\newblock {Asymptotic growth of multiplicity functions}.
\newblock {\em Journal of Pure and Applied Algebra}, 219:1045--1054, 2015.

\bibitem[CLO13]{cox}
D.~Cox, J.~Little, and D.~O'Shea.
\newblock {\em Ideals, varieties, and algorithms: an introduction to
  computational algebraic geometry and commutative algebra}.
\newblock Springer Science \& Business Media, 2013.

\bibitem[CN08]{nim}
G.~Colom\'e-Nin.
\newblock {\em Multigraded Structures and the Depth of Blow-up Algebras}.
\newblock PhD thesis, Universitat de Barcelona, 2008.

\bibitem[CSS19]{CSS}
S.D. Cutkosky, P.~Sarkar, and H.~Srinivasan.
\newblock Mixed multiplicities of filtrations.
\newblock {\em Transactions of the American Mathematical Society},
  372(9):6183--6211, 2019.

\bibitem[Cut13]{DC2}
S.~D. Cutkosky.
\newblock Multiplicities associated to graded families of ideals.
\newblock {\em Algebra \& Number Theory}, 7(9):2059--2083, 2013.

\bibitem[{Das}22]{das2}
S.~{Das}.
\newblock {Epsilon multiplicity for Noetherian graded algebras}.
\newblock {\em Illinois Journal of Mathematics}, 66(2):245--269, 2022.

\bibitem[Eis13]{eisenbud}
D.~Eisenbud.
\newblock {\em Commutative algebra: with a view toward algebraic geometry},
  volume 150.
\newblock Springer Science \& Business Media, 2013.

\bibitem[Gho16]{ghosh}
D.~Ghosh.
\newblock Asymptotic linear bounds of {C}astelnuovo--{M}umford regularity in
  multigraded modules.
\newblock {\em Journal of Algebra}, 445:103--114, 2016.

\bibitem[HH11]{HH}
J.~Herzog and T.~Hibi.
\newblock {\em Monomial ideals}.
\newblock Springer, 2011.

\bibitem[HHRT97]{HHRT}
M.~Herrmann, E.~Hyry, J.~Ribbe, and Z.~Tang.
\newblock Reduction numbers and multiplicities of multigraded structures.
\newblock {\em Journal of Algebra}, 197(2):311--341, 1997.

\bibitem[HHT07]{HHT}
J.~Herzog, T.~Hibi, and N.~V. Trung.
\newblock Symbolic powers of monomial ideals and vertex cover algebras.
\newblock {\em Advances in Mathematics}, 210(1):304--322, 2007.

\bibitem[HPV08]{J}
J.~Herzog, T.~J. Puthenpurakal, and J.~K. Verma.
\newblock Hilbert polynomials and powers of ideals.
\newblock {\em Mathematical Proceedings of the Cambridge Philosophical
  Society}, 145(3):623--642, 2008.

\bibitem[HS06]{HS}
C.~Huneke and I.~Swanson.
\newblock {\em Integral closure of ideals, rings, and modules}, volume~13.
\newblock Cambridge University Press, 2006.

\bibitem[JMn13]{JM}
J.~Jeffries and J.~Monta\~{n}o.
\newblock The $j$-multiplicity of monomial ideals.
\newblock {\em Mathematical Research Letters}, 20:729--744, 12 2013.

\bibitem[KV10]{KV}
D.~Katz and J.~Validashti.
\newblock Multiplicities and {R}ees valuations.
\newblock {\em Collectanea mathematica}, 61(1):1--24, 2010.

\bibitem[{Woo}15]{KW}
K.~{Woods}.
\newblock {Presburger arithmetic, rational generating functions and
  quasi-polynomials}.
\newblock {\em Journal of Symbolic Logic}, 80(2):433--449, 2015.

\end{thebibliography}

\end{document}